\documentclass[11pt]{article}
\usepackage{amsmath, amssymb, amscd, amsthm}
\title{The Area of Polynomial Images and Preimages}
\date{\today}
\author{{\Large{Edward Crane}} \thanks{Supported by an EPSRC research studentship}\vspace{6pt}\\ Department of Pure Mathematics and Statistics, University of Cambridge,\\ Cambridge CB3 0WA, United Kingdom.}

\newtheorem{theorem}{Theorem}
\newtheorem{lemma}{Lemma}

\newtheorem{definition}{Definition}

\newcommand{\ca}{\textup{cap}}
\begin{document}
\maketitle

\begin{abstract}Let $p$ be a monic complex polynomial of degree $n$ and let $K$ a measurable subset of the complex plane. We show that the area of $p(K)$, counted with multiplicity, is at least $\pi n \left( \textup{Area}(K)/\pi\right)^n$ and that \[\frac{\textup{Area}\left(p^{-1}(K)\right)}{\pi} \,\le\, \left(\frac{\textup{Area}(K)}{\pi}\right)^{1/n}\,.\] Both bounds are sharp. The special case of the latter result in which $K$ a disc was proved by P\'olya in 1928. We use Carleman's isoperimetric inequality relating the conductance and area for plane condensers. We include a summary of the necessary potential theory.
\end{abstract}
\maketitle
\section{Introduction and Statement of Results}
When $p$ is a polynomial of degree $n$ over $\mathbb{C}$, the set \[E(p,r) = \{z \in \mathbb{C}: |p(z)| = r^n \}\]is called a lemniscate, after Bernoulli's lemniscate $\{z \in \mathbb{C}: |z^2 - 1| = 1\}$. It is natural to ask how large the set enclosed by a lemniscate can be. We need to make some normalisation for this to be meaningful; the simplest is to ask for $p$ to be monic. Perhaps surprisingly, the area of the lemniscate $E(p,1)$ is then bounded independently of $p$. In fact the following sharp inequality was proved by P\'olya in 1928, \cite{Pol}.

\begin{theorem}[P\'olya's inequality]\label{T: disc}{\quad}\\
Let $p$ be a monic polynomial of degree $n$ over $\mathbb{C}$ and let $D$ be a disc in $\mathbb{C}$. Then the Euclidean area of $p^{-1}(D)$ is at most $\pi \left(\frac{\textup{Area}(D)}{\pi}\right)^{1/n}$, with equality only when $p: z \mapsto a (z-b)^n + c$ and the center of $D$ is $c$, the unique critical value of $p$.  
\end{theorem}

It is natural to ask whether we can get a larger preimage by fixing the area of $D$ but allowing its shape to vary. The main theorem of this paper is that we cannot.

\begin{theorem} \label{T: main} Let $p$ be a monic polynomial of degree $n$ over $\mathbb{C}$. Let $K$ be any measurable subset of the plane. Then \[\textup{Area}(p^{-1}(K)) \le \pi \left(\frac{\textup{Area}(K)}{\pi} \right)^{1/n}\,,\] with equality if and only if $K$ is (up to sets of measure zero) a disc and $p$ has a unique critical value at the center of that disc. 
\end{theorem}

Theorem \ref{T: main} is a simple consequence of the following stronger theorem:

\begin{theorem}\label{T: multiplicity}
Let $p$ be a monic polynomial of degree $n$ over $\mathbb{C}$, and $K$ be any measurable subset of the plane. Define the multiplicity $n(z, p, K)$ to be the number of $p$--preimages of $z$ in $K$, counted according to their valency. Then the area of $p(K)$ counted with multiplicity satisfies \[ \int_{\mathbb{C}} n(z, p, K) \, d\textup{A}(z) \; = \; \int_K |p'(z)|^2 \,d\textup{A} \;\ge\; n \pi \left( \frac{\textup{Area}(K)}{\pi}\right)^n \, ,\] with equality if and only if $K$ is (up to sets of measure zero) a disc and $p$ has a unique critical value at the center of that disc.
\end{theorem}

For a compact set $K \subset \mathbb{C}$, define \[\rho(K) = \frac{\textup{Area}(K)}{\pi\, \ca(K)^2}\,,\] where $\ca(K)$ is the logarithmic capacity of $K$. Then $\rho(K)$ is a measure of the \emph{roundness} of $K$: we will see that $\rho(K) \in [0,1]$, and $\rho(K)=1$ if and only if $K$ is a full--measure subset of a disc. We use $\rho$ to formulate the following scale-invariant version of theorem \ref{T: main}.

\begin{theorem}\label{T: roundness}
If $p$ is any complex polynomial of degree $n$, not necessarily monic, and $K$ is any compact subset of the plane, then \[\rho(p^{-1}(K)) \le \rho(K)^{1/n}\,.\]
\end{theorem}

This corollary is sharp for each value of $\rho$. To see this we can take $K$ to be the union of the unit disc with a radial line segment, and $p: z \mapsto z ^n$, so that we still get equality in theorem \ref{T: main}. 

In section \ref{S: Capacity} we give a quick introduction to the potential theory that we will need. In section \ref{S: Isoperimetric} we discuss some isoperimetric inequalities and their relationship with P\'olya's inequality (theorem \ref{T: disc}), which we prove since it is an important ingredient in the proof of theorem \ref{T: main}. Theorems \ref{T: main}, \ref{T: multiplicity} and \ref{T: roundness} are proved in section \ref{S: Proof}.

A survey of area estimates for lemniscates has recently been given by Lubinsky \cite{Lub}, with a view towards applications in the convergence theory of Pad\'e approximation. In \cite{EH},  Eremenko and Hayman made progress on the related problem of bounding the length of $E(p,r)$. Fryntov and Rossi \cite{FR} have obtained the sharp analogue of P\'olya's inequality (theorem \ref{T: disc}) bounding the hyperbolic area of the preimage of a hyperbolic disc under a finite Blaschke product. This raises the question of finding the sharp Blaschke product analogues of theorems \ref{T: main} and \ref{T: multiplicity}. 

 The author thanks Assaf Naor and Ben Green for posing the question that led to this paper, and his PhD supervisor Keith Carne for useful conversations.

\section{Capacity of plane subsets and condensers}\label{S: Capacity}

\begin{definition}: A \emph{plane condenser} is a pair $(E , B)$ of subsets of $\mathbb{C}$, where $E \neq \mathbb{C}$ is open and $B$ is a non-empty closed subset of $E$.
\end{definition}

The terminology arises from the fact that a pair of conducting cylinders with cross-section $\partial E$ and $B$ respectively could be used as a condenser (or capacitor). The \emph{capacity} of the condenser $(E, B)$ is physically the capacitance per unit length of an infinitely long pair of such cylinders. The same quantity describes the \emph{conductance} between $\partial E$ and $B$ of an isotropic resistor consisting of a plate in the shape of $E \setminus B$. We compute this conductance by considering the electrical potential $f$ that would be induced in $U$ if we were to connect $\partial E$ to an electrical potential $0$ and $B$ to potential $+1$. Given $f$, we can compute the current $-\nabla f$ that would flow in response to the potential $f$, and the power consumed is proportional to \[L(f) = \int_{E \setminus B} |\nabla f|^2 \,d\text{A}\,,\] We expect a physical potential $f$ to minimise $L(f)$ over all possible potential functions satisfying the given boundary conditions. Calculus of variations tells us that if there is an extremal $f$, it must be harmonic on $E \setminus B$. We call such an $f$ a Green's function for the condenser. A Green's function only exists if the boundary $\partial E \cup \partial B$ is regular for the Dirichlet problem, but we avoid this difficulty by defining \[\ca(E,B)\; =\; \frac{1}{4\pi}\; \inf L(f)\,,\] where the infimum is taken over all continuously differentiable $f: \mathbb{C} \to \mathbb{R}$ such that $f = 0$ on $\mathbb{C} \setminus E$ and $f = 1$ on $B$. We call such functions \emph{admissible} for the condenser $(E,B)$. Note that $\ca(E,B)$ may be zero, as it is when $B$ is a finite set. From the definition it is immediate that capacity is monotonic, i.e. \[E \subseteq F \text{ and } B \supseteq C \quad \implies \quad \ca(E,B) \ge \ca(F,C). \] 
\begin{lemma}\label{L:pullback} Suppose that on some open set $U \subset \mathbb{C}$ we have an analytic function $\psi$ such that each point of $E$ has exactly $n$ preimages in $U$, counted according to valency.
 Then  \[\ca (\psi^{-1}(E), \psi^{-1}(B))= n \; \ca (E,B)\,.\]
\end{lemma}
\begin{proof} Suppose that $f: \mathbb{C} \to \mathbb{R}$ is any admissible function for $(E,B)$ . The hypothesis implies that the restriction of $\psi$ to $\psi^{-1}(E)$ is a proper map, so we can extend $f \circ \psi$ to get an admissible function for $(\psi^{-1}(E), \psi^{-1}(B))$ by giving it the value $0$ outside $U$. Since $\psi$ is almost everywhere conformal, \begin{eqnarray*} L(f \circ \psi) &=& \int_{\psi^{-1}(E \setminus B)} |\nabla(f \circ \psi)(z)|^2 \,d\textup{Area}(z)\\ &=& \int_{\psi^{-1}(E \setminus B)} |(\nabla f )(\psi(z))|.|\psi'(z)|^2 \, d\textup{Area}(z)\\ &=& n \int_{E \setminus B} |\nabla f (w)|^2 d\textup{Area}(w) \quad = \quad n\,L(f)\,.\end{eqnarray*}
\end{proof}
In particular the capacity is a \emph{conformal invariant} of condensers: if $\varphi: E \to \mathbb{C}$ is a univalent function then \[\ca(E,B) = \ca(\varphi(E), \varphi(B))\,.\]  For example, if $(E \setminus B)$ is a ring domain then its modulus is $1/(4 \pi \ca(E,B)$. 

Let $K$ be a compact set in the plane. A Green's function for $K$ is a continuous function $f: \mathbb{C} \to \mathbb{R}$, zero on $K$ and harmonic on $\mathbb{C} \setminus K$, with $f(z) = \log |z| -\log t + o(1)$ as $z \to \infty$. If $K$ has a Green's function then the logarithmic capacity of $K$ is defined to be $\ca(K) = t$. for general $K$ we define $\ca(K) = \inf \ca(J)$ over all compact sets $J \supset K$ with regular boundary for the Dirichlet problem on $\mathbb{C} \setminus J$. By pulling back Green's functions, it is easy to verify that if $p$ is a monic polynomial of degree $n$ then \[ \ca ( p^{-1}(B)) = \ca(B)^{1/n}\,.\]

\section{Isoperimetric Inequalities}\label{S: Isoperimetric}

A relationship between capacity and 2-dimensional Lebsegue measure is given by the following `isoperimetric' inequality:

\begin{theorem}\label{T: isoperimetric inequality for condensers}\textup{(Carleman, 1918)}\\
\[\frac{1}{\ca(E,B)} \le \log \left(\frac{\textup{Area}(E)}{\textup{Area}(B)}\right)\,,\] with equality iff $E$ and $B$ are concentric discs.
\end{theorem}

The proof of Carleman's inequality uses the fact that the Dirichlet integral $L(f)$ does not increase when $f$ is replaced by its Schwarz symmetrization, the function $S(f)$ whose superlevel sets are concentric discs with the same area as the corresponding level sets of $f$. For details, see the classic book of P\'olya and Szeg\"{o}, \cite{PS}, or \cite{Ban} for a more modern account.

Taking $E = B(0,R)$ and then letting $R \to \infty$ in Carleman's isoperimetric inequality yields the following well-known isoperimetric theorem for logarithmic capacity. For a simple proof, including the equality case, see theorem 5.3.5 in \cite{Ran}.

\begin{theorem}\label{T: isoperimetric inequality for logarithmic capacity}
For any compact set $K \subset \mathbb{C}$, \[\textup{Area}(K)\; \le\; \pi \,\ca(K)^2\,,\] with equality if and only if $K$ is a disc. 
\end{theorem}

We have now collected everything we need to prove P\'olya's inequality, theorem \ref{T: disc}. The capacity of the disc $D$ is precisely the radius of $D$, so \[\ca (D) = \left(\frac{\textup{Area}(D)}{\pi}\right)^{1/2}\,,\] \[\ca(p^{-1}(D)) = \left(\frac{\textup{Area}(D)}{\pi}\right)^{1/2n}\,,\] and, applying theorem \ref{T: isoperimetric inequality for logarithmic capacity}, \[\textup{Area}(p^{-1}(D)) \le \pi \left(\frac{\textup{Area}(D)}{\pi}\right)^{1/n}\,,\] as required. In view of the strong link between logarithmic capacity and polynomials, theorems \ref{T: disc} and \ref{T: isoperimetric inequality for logarithmic capacity} are virtually equivalent. In \cite{Lub}, P\'olya's inequality is proved using Gronwall's area formula, and used to deduce the isoperimetric inequality for logarithmic capacity.

\section{Proof of theorems 2 and 3}\label{S: Proof}

\begin{lemma}\label{L: integrated Carleman} For any complex polynomial $g$ of degree $d$, \[\int_{\mathbb{C}} |g(w)|\, {\bf{1}}_{|g(w)| \le x}\; d\textup{A} \; \ge \; \frac{2x}{d+2} \; \textup{Area}(\{w \in \mathbb{C} : |g(w) \le x\})\,.\]
\end{lemma}
\begin{proof}
By lemma \ref{L:pullback}, we have \[\ca\left( g^{-1}(B(0,x)),g^{-1}(B(0,s))\right) = \frac{d}{2 (\log x - \log s)}\,.\] Theorem \ref{T: isoperimetric inequality for condensers} gives \[\frac{\textup{Area}\left(\{w \in \mathbb{C} : s \le |g(w)| \le x\right\})}{\textup{Area}\left(\{w \in \mathbb{C} : |g(w) \le x\}\right)} \; \ge \; 1 - \left(\frac{s}{x}\right)^{2/d} \,,\] so
\begin{eqnarray*}
\int_{\mathbb{C}} |g(w)| \,{\bf{1}}_{|g(w)| \le x}\; d\textup{A} & = & \int_{0}^{x} \textup{Area}\left(\{w \in \mathbb{C} : s \le |g(w)| \le x\right\}) \, ds\\ & \ge &  \textup{Area}\left(\{w \in \mathbb{C} : |g(w) \le x\}\right) \int_{0}^{x} 1 -  \left(\frac{s}{x}\right)^{2/d} \,ds \\& = & \frac{2x}{d+2} \; \textup{Area}(\{w \in \mathbb{C} : |g(w) \le x\})\,.
\end{eqnarray*}
\end{proof}

 Now fix a monic polynomial $p$ and $A >0$. Among all measurable sets $K$ with $\textup{Area}(K) = A$, the Dirichlet integral \[\int_K |p'(w)|^2 d\textup{Area}(w)\] is minimised when $K$ is the sublevel set \[K_t = \{w \in \mathbb{C}: |p'(w)|^2 \le t\}.\] Here $t$ is determined uniquely by the condition $\textup{Area}(K_t)  = A$. The polynomial $z \mapsto (p'(z) /n)^2$ is monic, with degree $2n-2$, so theorem \ref{T: disc} gives \[A = \textup{Area}(K_t) \le \pi \, \left(\frac{\pi (t/n^2)^2}{\pi}\right)^{1/(2n-2)}\,.\] Rearranging this we have \[t \ge n^2 \left(\frac{A}{\pi}\right)^{n-1}\,.\]  
Now we apply lemma \ref{L: integrated Carleman} to the polynomial $g = (p')^2$ to obtain \begin{eqnarray*} \int_{K_t} |p'(w)|^2 d\textup{A}(w) & = &\int_\mathbb{C} |p'(w)|^2 {\bf{1}}_{|p'(w)|^2 \le t} \, d\textup{A}(z) \\
& \ge &\frac{2t}{2n} \textup{Area}(K_t) = \frac{tA}{n}\\ &\ge& n \pi \left(\frac{A}{\pi}\right)^n\,.
\end{eqnarray*}
For equality, we must have equality in our application of P\'olya's inequality, so $p$ must be $p: z \mapsto (z-b)^n + c$, and $K$ can differ from disc $K_t$ at most by a set of 2--dimensional Lebesgue measure zero. This completes the proof of theorem \ref{T: multiplicity}.

To obtain theorem \ref{T: main}, observe that a monic polynomial $p$ maps $p^{-1}(K)$ onto $K$ with multiplicity $n$ everywhere, so \[\textup{Area}(K) = \frac{1}{n} \int_{p^{-1}(K)} |p'(w)|^2\, d\textup{A}(w) \,.\]
Finally, theorem \ref{T: roundness} is obtained by dividing both sides of the inequality of theorem \ref{T: main} by $\ca(K)^2 = \ca(p^{-1}(K))^2$.

\end{document}